\documentclass{amsart}
\usepackage{amsmath,amssymb,amsfonts,graphicx, hyperref}

\usepackage{array,rotating}

\title{The contact polytope of the Leech lattice (complete version)}

\author{Mathieu Dutour Sikiri\'c}
\address{M.~Dutour Sikiri\'c, Rudjer Boskovi\'c Institute, Bijenicka 54, 10000 Zagreb, Croatia}
\email{mdsikir@irb.hr}

\author{Achill Sch\"urmann}
\address{A.~Sch\"urmann, Delft Institute of Applied Mathematics,
  Technical University of Delft, P.O. Box 5031, 2600 GA Delft, The
  Netherlands}
\email{a.schurmann@tudelft.nl}

\author{Frank Vallentin} 
\address{F.~Vallentin, Delft Institute of Applied Mathematics,
  Technical University of Delft, P.O. Box 5031, 2600 GA Delft, The
  Netherlands}
\email{f.vallentin@tudelft.nl}



\date{February 16, 2010}

\newcommand{\RR}{\ensuremath{\mathbb{R}}}
\newcommand{\NN}{\ensuremath{\mathbb{N}}}

\newcommand{\ZZ}{\ensuremath{\mathbb{Z}}}

\newtheorem{theorem}{Theorem}

\DeclareMathOperator{\conv}{conv}

\DeclareMathOperator{\Min}{Min}
\DeclareMathOperator{\Cont}{C}
\DeclareMathOperator{\Co}{Co}
\DeclareMathOperator{\DV}{V}

\begin{document}

\begin{abstract}
  The contact polytope of a lattice is the convex hull of its shortest
  vectors.  In this paper we classify the facets of the contact
  polytope of the Leech lattice up to symmetry. There are
  $1,197,362,269,604,214,277,200$ many facets in $232$ orbits.
\end{abstract}

\maketitle

\section{Introduction}

An $n$-dimensional \textit{lattice} $L$ is a discrete subgroup of the $n$-dimensional Euclidean space $\RR^n$ of the form $L = \{\sum_{i = 1}^n \alpha_i b_i : \alpha_1, \ldots, \alpha_n \in \ZZ\}$ where $b_1, \ldots, b_n$ is a basis of $\RR^n$.  By $\lambda(L)$ we denote the Euclidean length of non-zero shortest vectors of $L$ and we denote by $\Min L$ the set of \textit{shortest vectors}.

Every lattice comes with two important polytopes: The \textit{contact polytope} of $L$ is the convex hull of its shortest vectors
\begin{equation*}
\Cont(L) = \conv\left\{v : v \in \Min L\right\},
\end{equation*}
and the \textit{Voronoi cell} of $L$ is the region of points which are closer to the origin than to other lattice points
\begin{equation*}
\DV(L) = \left\{x\in \RR^n : x \cdot v \leq \frac{1}{2} v \cdot v \text{ for all }v\in L\right\}.
\end{equation*}

Maybe one of the most remarkable lattices is the $24$-dimensional Leech lattice $\Lambda_{24}$. It has $196560$ shortest vectors which is the highest possible number in dimension~$24$. Its \textit{orthogonal group}, i.e.\ the group of orthogonal transformations preserving the lattice is the Conway group~$\Co_0$. It has $2^{22} \cdot 3^9 \cdot 5^4 \cdot 7^2 \cdot 11 \cdot 13 \cdot 23 = 8,315,553,613,086,720,000$ elements and it is connected to many sporadic simple groups. We refer to the book \cite{CS} by Conway and Sloane for an extensive treatment of the Leech lattice.

Borcherds, Conway, Parker, Queen, Sloane \cite[Chapter 23, Chapter 25]{CS} determine the vertices of the Voronoi cell of the Leech lattice. The Voronoi cell tiles the space~$\RR^n$ by translations; this gives the \textit{Voronoi cell tiling} of $\RR^n$.  So, in the context of the Voronoi cell it is natural to consider orbits under the \textit{isometry group} (the group generated by the orthogonal group of the Leech lattice together with lattice translations) acting on the Voronoi cell tiling. We denote the isometry group of the Leech lattice by $\Co_{\infty}$. There are $307$ orbits of vertices in the Voronoi cell tiling under the action of~$\Co_{\infty}$.

In this paper we determine the facets and their incidence relations of the contact polytope of the Leech lattice. We get the following result.

\begin{theorem}
  There are $232$ orbits of facets of $\Cont(\Lambda_{24})$ under
  $\Co_0$.
\end{theorem}

The contact polytope and the Voronoi cell are related. To see this relation, we consider \begin{equation*}
  \Cont(L)^* = \left\{x \in \RR^n : x \cdot v \leq \frac{1}{2} v \cdot v \text{ for all }v\in \Min L\right\},
\end{equation*}
which is the standard polar polytope scaled by a factor of $\frac{1}{2} \lambda(L)^2$. The faces of $\Cont(L)$ and of $\Cont(L)^*$ are in bijection. The bijection reverses the inclusion relation: $k$-dimensional faces of $\Cont(L)$ correspond to $(n-k)$-dimensional faces of $\Cont(L)^*$. In particular, vertices of $\Cont(L)^*$ correspond to facets of $\Cont(L)$. For these notions we refer to the standard literature on polytope theory, e.g.\ the book by Ziegler~\cite{ziegler-1995}.

We chose the scaling in the definition of $\Cont(L)^*$ so that it contains $\DV(L)$. In the case of the Leech lattice some vertices of $\DV(\Lambda_{24})$ and $\Cont(\Lambda_{24})^*$ are shared. As a side remark: One has equality $\Cont(L)^* = \DV(L)$ if and only if $L$ is a root lattice, see Rajan, Shende \cite{RS}.

\begin{theorem}
\label{th:thm2}
  $164$ orbits of vertices of $\Cont(\Lambda_{24})^*$ are also orbits of
  vertices of $\DV(\Lambda_{24})$. They are listed in Table~1 in the complete version of the paper \cite{table}.
  The additional $68$ orbits of vertices are listed in Table~2 of \cite{table}.
\end{theorem}

We classify the shared vertices in Section~\ref{sec:shared} and give them in Table~1 of \cite{table}. In Section~\ref{sec:additional} we classify the additional vertices of $\Cont(\Lambda_{24})^*$ which are not vertices of $\DV(\Lambda_{24})$. We conclude the paper by Section~\ref{sec:techniques} where we briefly explain our computational techniques.

The data presented here is also electronically available from \cite{contact}.

\section{Shared vertices}
\label{sec:shared}

In this section we explain the notation used in Table~1 of \cite{table} which contains the $164$ orbits of shared vertices mentioned in Theorem~\ref{th:thm2}.

The vertices of the Voronoi cell of a lattice are centers of \textit{empty spheres}, i.e.\ spheres $S(x,\|x\|)$ with center $x$ and radius $\|x\|$ which contain lattice points on their boundary but not in their interior. The convex hull of lattice points on the boundary of such an empty sphere is called the \textit{Delone cell} of the vertex $x$.

The Delone cells of the Leech lattice are classified by Borcherds, Conway, Parker, Queen, Sloane \cite[Chapter 23, Chapter 25]{CS} up to the action of the isometry group $\Co_{\infty}$. For this classification they use Coxeter-Dynkin diagrams.

A \textit{Coxeter-Dynkin diagram} with vertex-set $\{1,\dots, N\}$ is a symmetric $N\times N$ matrix $(m_{ij})_{1\leq i,j\leq N}$ with ones on the diagonal and $m_{ij} \geq 2$ if $i\not= j$ and $m_{ij}\in \NN\cup \{\infty\}$.

A Coxeter-Dynkin diagram is called \textit{simply laced} if $m_{ij}=2$, $3$ or $\infty$.  The \textit{Cartan matrix} of a Coxeter-Dynkin diagram $(m_{ij})_{1\leq i,j\leq N}$ is the matrix $M=\bigl(-\cos\,\frac{\pi}{m_{ij}}\bigr)_{1\leq i,j\leq N}$.  A Coxeter-Dynkin diagram is called \textit{spherical} if its Cartan matrix is positive definite and \textit{affine} if its Cartan matrix is positive semidefinite.  A Coxeter-Dynkin diagram is called \textit{decomposable} if we can partition its vertex-set into $S_1\cup S_2$ with $m_{ij}=2$ if $i\in S_1$ and $j\in S_2$. It is called \textit{indecomposable} otherwise. A Coxeter-Dynkin diagram~$D$ admits a unique decomposition into indecomposable Coxeter-Dynkin diagrams $D_1$, \dots, $D_r$, which we write as $D=D_1D_2\dots D_r$.  The classification of spherical and affine Coxeter-Dynkin diagrams is presented, for example, in Humphreys \cite[Section 2.4, 4.7]{humphreyscoxeter}. Here the famous $A-D-E$ diagrams show up, explained e.g.\ by Hazewinkel, Hesselink, Siersma, Veldkamp~\cite{HHSV}. The spherical, simply laced, indecomposable Coxeter-Dynkin diagrams are $a_n$ for $n\geq 1$, $d_n$ for $n\geq 4$ and $e_n$ for $6\leq n\leq 8$. Each diagram corresponds to an indecomposable affine diagram: $A_n$, $D_n$ and $E_n$. All these diagrams are pictured e.g.\ in \cite[Figure 23.1]{CS}.

In the Leech lattice, a Coxeter-Dynkin diagram $(m_{ij})_{1\leq  i,j\leq N}$ can be associated with a Delone cell with vertex-set $\{v_1, \dots, v_N\}$ by
\begin{equation*}
m_{ij}=\left\{\begin{array}{rcl}
1 &\text{if} &\|v_i -v_j\|^2=0,\\
2 &\text{if} &\|v_i -v_j\|^2=4,\\
3 &\text{if} &\|v_i -v_j\|^2=6,\\
\infty &\text{if} &\|v_i -v_j\|^2=8.
\end{array}\right.
\end{equation*}
As can been seen in Table~1 of \cite{table}, different Delone cells may have the same Coxeter-Dynkin diagram.

In Table~1 of \cite{table} the rows are sorted first by the squared length $\|v\|^2$ (third column) of the vertex $v$.  Second they are sorted by the size of the stabilizer of $v$ within the orthogonal group of the Leech lattice (fifth column), and then by the number of incident facets of $\Cont(\Lambda_{24})^*$ (fourth column).

In the second column we give the Coxeter-Dynkin diagrams of the associated Delone cell of $v$. Note that the diagrams are affine if and only if the squared length of $v$ equals $2$, the maximum among shared vertices. In all other cases they are spherical. Furthermore, in the spherical cases the number of incident facets is always equal to the minimum possible number of $24$. These observations follow from \cite[Chapter 23, Chapter 25]{CS}.

In the last column we give the MOG (Miracle Octad Generator) coordinates of representatives of each orbit which one has to multiply with $\alpha$ (sixth column). The MOG coordinates form a standard coordinate system for the Leech lattice. They are explained in~\cite[Chapter 11]{CS}.

There are $307$ orbits of vertices in the Voronoi cell tiling under the action of the isometry group $\Co_{\infty}$ of the Leech lattice. Our computation shows that there are $5297$ orbits of vertices of the single Voronoi cell $\DV(\Lambda_{24})$ under the action of the smaller, finite orthogonal group of the Leech lattice; $164$ of them are shared with $\Cont(\Lambda_{24})^*$.

\section{Additional vertices}
\label{sec:additional}

There are $68$ additional orbits of vertices of $\Cont(\Lambda_{24})^*$ which are not vertices of the Voronoi cell of the Leech lattice. These additional vertices are characterized by the fact that the distance to a closest lattice point is strictly less than the distance $\|v\|$ to the origin.

Table~2 of \cite{table} describes these $68$ orbits. Like in Table~1 of \cite{table} the rows are sorted (in this order) by the squared length $\|v\|^2$ (third column), the size of the stabilizer of $v$ within the orthogonal group of the Leech lattice (fifth column), and then by the number of incident facets (fourth column).

In the second column we give names for diagrams. The first row corresponds to an exceptional vertex which we explain below. The other $67$ rows correspond to graphs which we define later in Section~\ref{ssec:other}.

\subsection{The exceptional vertex}
\label{ssec:exceptional}

The first orbit of vertices is exceptional: Its squared norm $8/3 = 2.666\dots$ is substantially bigger than the squared norm of all other vertices which lie in the interval $[1.92, 2.25]$. Its incidence number of $552$ as well as the size of its stabilizer, which is the Conway group $\Co_3$, are also substantially larger than the values for the other vertices. This orbit of vertices is a scaled copy of the vectors of $\Lambda_{24}$, having Euclidean norm $\sqrt{6}$.

In the contact polytope $\Cont(\Lambda_{24})$ this exceptional vertex corresponds to a facet. Since it has maximum norm among all vertices the corresponding facet is closest to the origin and has the largest possible circumsphere among all other facets of $\Cont(\Lambda_{24})$. This solves a conjecture of Ballinger, Blekherman, Cohn, Giansiracusa, Kelly, Sch\"urmann \cite[Section 3.7]{BBCGKS}. We note that a similar calculation as the one presented here, solves the corresponding conjecture about the contact polytope of the $23$-dimensional lattice $O_{23}$, the shorter Leech lattice, which has $4600$ vertices.

The $23$-dimensional point configuration, given by the $552$ shortest vectors of the Leech lattice defining facets incident to the exceptional vertex, appears in several different contexts: It is universally optimal (Cohn, Kumar \cite{CK}), it defines $276$ equiangular lines (Lemmens, Seidel \cite{LS}), and it defines an extreme Delone cell (Deza, Laurent \cite[Chapter 16.3]{DL}). Moreover, it contains a wealth of remarkable substructures (see Cohn et. al. \cite{BBCGKS}), e.g.\ the highly-symmetric point configurations discussed in the next section, but also others, e.g.\ the one defined by the McLaughlin graph.

\subsection{The other vertices}
\label{ssec:other}

To the remaining $67$ orbits of vertices we associate a diagram as follows. Let $v$ be one of these vertices and let $w_1, \ldots, w_N$ be shortest vectors of the Leech lattice defining facets incident to $v$. Only the two inner products $1$ and $2$ occur between distinct vectors $w_i$ and $w_j$. So we can define a graph with vertex-set $\{1, \ldots, N\}$ and edge-set $\{\{i,j\} : w_i \cdot w_j = 1\}$; the other inner product $2$ defines non-edges.

Here again the graphs decompose into connected components where several of these occurring components are highly-symmetric and have been studied in other contexts. We discuss them below, the graphs $a_n$, $d_n$ and $e_n$ are already described in the previous section, and the remaining ones are in Figure~1.

The \textit{Higman-Sims graph} $\mathrm{HS}_{100}$ is the unique strongly regular graph with parameters $(100,22,0,6)$. See Brouwer, Cohen, Neumaier \cite[Chapter 13.1]{BCN}.

The \textit{Hoffman-Singleton graph} $\mathrm{HS}_{50}$ is the unique strongly regular graph with parameters $(50,7,0,1)$. See \cite[Chapter 13.1]{BCN}.

For the \textit{Johnson graph $J(7,4)$} see \cite[Chapter 9.1]{BCN}.

A \textit{$(k,g)$-cage} is a regular graph of valency $k$, girth $g$, which attains the minimum possible number of vertices. The $(5,6)$-cage (incidence graph of a projective plane $\mathrm{PG}(2,4)$) and the $(3,8)$-cage (\textit{Tutte-Coxeter graph}) are unique. See \cite[Chapter 6.9]{BCN} and Tutte \cite{T}.

The \textit{Coxeter graph} $\mathrm{Cox}$ is the unique distance regular graph with intersection array $\{3,2,2,1;1,1,1,2\}$. See \cite[Chapter 12.3]{BCN}.

In Figure~1 we list the remaining graphs. The vertices of these graphs only have degree one (white circles), degree two (sitting on edges, which are not depicted, but see below), or degree three (black circles).  We have three kinds of trees: $T^a_bc$ having $a+b+c+4$ vertices, $T^a_bc^d_e$ having $a+b+c+d+e+6$ vertices, and $T^a_bc^de^f_g$ having $a+b+c+d+e+f+g+8$ vertices; we have $12$ other graphs $G_{n,m}$ with $n$ vertices and $m$ edges. In Figure~1 the numbers on the edges show how many vertices of degree~2 sit on them, but in the following four cases we did not put these numbers: The graph $G_{24,30}$ has one vertex of degree~2 on every edge, $G_{25,30}$ is the Petersen graph which has one vertex of degree~2 on every edge, $G_{22,22}$ has three vertices of degree~2 on every edge, and the graph $G_{24,27}$ is the complete bipartite graph $K_{3,3}$ which has two vertices of degree~2 on every edge.

\begin{table}[htb]
\begin{tabular}{ccc}
\includegraphics[height = 2.5cm]{tabc.eps} & 
\includegraphics[height = 2.5cm]{tabcde.eps} & 
\includegraphics[height = 2.5cm]{tabcdefg.eps} \\
$T^a_bc$ & $T^a_bc^d_e$ & $T^a_bc^de^f_g$ \\[0.3cm]
\includegraphics[height = 2.4cm]{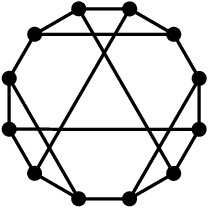} & 
\includegraphics[height = 3cm]{g2124.eps} & 
\includegraphics[height = 2.8cm]{g1920.eps} \\
$G_{24,30}$ & $G_{21,24}$ & $G_{19,20}$ \\[0.3cm]
\includegraphics[height = 3cm]{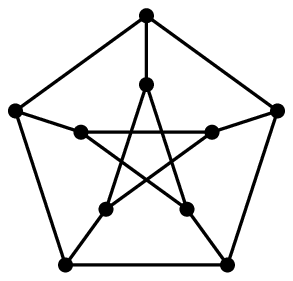} & 
\raisebox{0.2cm}{\includegraphics[height = 2.8cm]{g1717.eps}} &
\raisebox{0.3cm}{\includegraphics[height = 1.9cm]{g1515.eps}}\\
$G_{25,30}$ & $G_{17,17}$ & $G_{15,15}$\\[0.3cm]
\raisebox{0.5cm}{\includegraphics[height = 1.4cm]{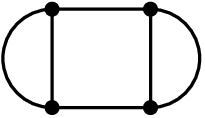}} & 
\includegraphics[height = 2.5cm]{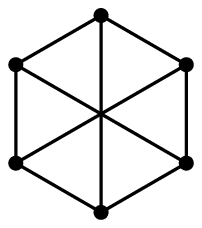} & 
\raisebox{0.5cm}{\includegraphics[height = 1.1cm]{g1414.eps}} \\
$G_{22,22}$ & $G_{24,27}$ & $G_{14,14}$ \\[0.3cm]
\includegraphics[height = 2.5cm]{g2020.eps} & 
\raisebox{0.3cm}{\includegraphics[height = 1.9cm]{g1818.eps}} & 
\raisebox{0.2cm}{\includegraphics[height = 2.8cm]{g2121.eps}} \\
$G_{20,20}$ & $G_{18,18}$ & $G_{21,21}$
\end{tabular}
\bigskip

\textsc{Figure 1. Diagrams}

\end{table}

\section{Computational techniques}
\label{sec:techniques}

Computing the vertices of $\Cont(\Lambda_{24})^*$ from its facets is called a \textit{polyhedral representation conversion problem}. A direct application of standard programs like Fukuda's {\tt cdd} \cite{cdd} or Avis' {\tt lrs} \cite{lrs} for this conversion is not feasible due to the extremely large number of vertices.

In order to exploit the symmetries of $\Cont(\Lambda_{24})^*$, we use the {\em adjacency decomposition method} which is surveyed in Bremner, Dutour Sikiri\'c, Sch\"urmann \cite{bremner}. An implementation by the first author is available from \cite{polyhedral}.

The adjacency decomposition method computes a complete list of inequivalent vertex representatives. First one computes an initial vertex by solving a linear program and inserts it into the list of orbit representatives.  From any such representative, we compute the list of adjacent vertices, and if they give a new orbit, we insert it into the list of representatives. After finitely many steps all orbits have been treated.  Computing adjacent vertices is a lower-dimensional representation conversion problem. So this method can be applied recursively.

For $\Cont(\Lambda_{24})^*$ we had to come up with two case-specific insights:

From \cite{BBCGKS} it is known that the exceptional vertex of Section~\ref{ssec:exceptional} is indeed a vertex of $\Cont(\Lambda_{24})^*$. We used it as starting vertex of the adjacency decomposition method.

For checking isomorphy and for computing stabilizers we used the following standard strategy: we characterize a vertex of $\Cont(\Lambda_{24})^*$ by the set of its incident facets and we represent the symmetry group $\Co_0$ as a permutation group acting on the $196560$ shortest vectors of the Leech lattice. Then, we use the backtracking algorithm by Leon \cite{leon-1991,leon-1997} implemented in \cite{GAP}. This worked reasonably fast for all the cases except for the two orbits of vertices having the same Coxeter-Dynkin diagram $a_1^{25}$. The stabilizer of the corresponding Delone cell under the isometry group $\Co_{\infty}$ is the Mathieu group $M_{24}$. Under the action of $M_{24}$ the $25$ vertices of the Delone cell split into two orbits of size $1$ and $24$. Hence, these two orbits correspond to two distinct orbits of vertices of $\Cont(\Lambda_{24})^*$, one having stabilizer $M_{24}$ and the other having stabilizer $M_{23}$. The backtracking algorithm of GAP could not decide in reasonable time whether or not two vertices with the same Coxeter-Dynkin diagram $a_1^{25}$ are in the same orbit. So we used the third method of Section~3.5 of \cite{DSV08} to resolve this problem.

\section*{Acknowledgements}

We started this research during the Junior Trimester Program (February 2008--April 2008) on ``Computational Mathematics' at the Hausdorff Institute of Mathematics (HIM) in Bonn. Then, part of this research was done at the Mathematisches Forschungsinstitut Oberwolfach during a stay within the Research in Pairs Programme from May 3, 2009 to May 16, 2009. We thank both institutes for their hospitality and support. The work of the first author has been supported by the Croatian Ministry of Science, Education and Sport under contract 098-0982705-2707. The second and the third author were supported by the Deutsche Forschungsgemeinschaft (DFG) under grant SCHU 1503/4-2.

\begin{sidewaystable}
\begin{center}

\begin{tabular}{c|c|c|>{$}c<{$}|>{$}c<{$}|c|>{$}r<{$}>{$}r<{$}>{$}r<{$}>{$}r<{$}|>{$}r<{$}>{$}r<{$}>{$}r<{$}>{$}r<{$}|>{$}r<{$}>{$}r<{$}>{$}r<{$}>{$}r<{$}}

& name & $\|v\|^2$ & $N$ & $g$ & $\alpha$ & \multicolumn{4}{c}{M} & \multicolumn{4}{c}{O} & \multicolumn{4}{c}{G}\\
\hline

1. & $A_1^{24}$ & $2$& 48 & 20891566080 &  $\frac{1}{2}$ &  0 &  2 &  0 &  2 &  0 &  -2 &  2 &  0 &  -2 &  0 &  0 &  2 \\ 
& & & & & &  -2 &  0 &  -2 &  0 &  2 &  0 &  0 &  2 &  -2 &  0 &  2 &  4 \\ 

2. & $D_{4}^6 $ & $2$& 30 & 5760 &  $\frac{1}{6}$ &  1 &  5 &  1 &  7 &  1 &  -9 &  7 &  -5 &  -5 &  1 &  1 &  5 \\ 
& & & & & &  -7 &  1 &  -3 &  3 &  5 &  1 &  3 &  5 &  -7 &  3 &  5 &  9 \\ 

3. & $A_{5}^4  D_{4} $ & $2$& 29 & 864 &  $\frac{1}{6}$ &  4 &  6 &  0 &  6 &  2 &  -10 &  2 &  -6 &  -2 &  8 &  0 &  -2 \\ 
& & & & & &  -4 &  0 &  -6 &  2 &  2 &  4 &  6 &  4 &  -2 &  2 &  6 &  10 \\ 

4. & $E_{6}^4$ & $2$& 28 & 36 &  $\frac{1}{12}$ &  2 &  8 &  6 &  8 &  2 &  -18 &  16 &  -12 &  -14 &  2 &  0 &  8 \\
\end{tabular}

\bigskip

\textsc{Table 1. Shared vertices} 
\end{center}
\end{sidewaystable}

\begin{sidewaystable}
\begin{center}

\begin{tabular}{c|c|c|>{$}c<{$}|>{$}c<{$}|c|>{$}r<{$}>{$}r<{$}>{$}r<{$}>{$}r<{$}|>{$}r<{$}>{$}r<{$}>{$}r<{$}>{$}r<{$}|>{$}r<{$}>{$}r<{$}>{$}r<{$}>{$}r<{$}}

& name & $\|v\|^2$ & $N$ & $g$ & $\alpha$ & \multicolumn{4}{c}{M} & \multicolumn{4}{c}{O} & \multicolumn{4}{c}{G}\\
\hline
21. & $D_{16}  E_{8} $ & $2$& 26 & 2 &  $\frac{1}{30}$ &  0 &  24 &  12 &  22 &  0 &  -50 &  36 &  -28 &  -28 &  12 &  2 &  20 \\ 
& & & & & &  -42 &  4 &  -16 &  12 &  26 &  6 &  20 &  18 &  -30 &  24 &  16 &  44 \\ 

22. & $D_{16}  E_{8} $ & $2$& 26 & 2 &  $\frac{1}{30}$ &  -2 &  20 &  10 &  20 &  -2 &  -44 &  38 &  -28 &  -26 &  10 &  2 &  14 \\ 
& & & & & &  -42 &  10 &  -14 &  10 &  26 &  4 &  22 &  16 &  -34 &  24 &  18 &  52 \\ 

23. & $D_{16}  E_{8} $ & $2$& 26 & 1 &  $\frac{1}{30}$ &  27 &  29 &  -3 &  27 &  3 &  -43 &  15 &  -27 &  -9 &  45 &  1 &  -13 \\ 
& & & & & &  -27 &  -5 &  -31 &  15 &  11 &  13 &  27 &  21 &  -11 &  11 &  31 &  49 \\ 

24. & $D_{16} E_{8} $ & $2$ & 26 & 1 &  $\frac{1}{30}$ &  25 &  29 &  -3 &  25 &  3 &  -41 &  11 &  -29 &  -9 &  47 &  1 &  -15 \\
\end{tabular}

\bigskip

\textsc{Table 1 (contd.). Shared vertices} 
\end{center}
\end{sidewaystable}

\begin{sidewaystable}
\begin{center}

\begin{tabular}{c|c|c|>{$}c<{$}|>{$}c<{$}|c|>{$}r<{$}>{$}r<{$}>{$}r<{$}>{$}r<{$}|>{$}r<{$}>{$}r<{$}>{$}r<{$}>{$}r<{$}|>{$}r<{$}>{$}r<{$}>{$}r<{$}>{$}r<{$}}

& name & $\|v\|^2$ & $N$ & $g$ & $\alpha$ & \multicolumn{4}{c}{M} & \multicolumn{4}{c}{O} & \multicolumn{4}{c}{G}\\
\hline
41. & $a_{1} a_{10} d_{14}$ & $\frac{3716}{1859}$ & 25 & 1 &  $\frac{1}{143}$ &  122 &  132 &  -12 &  122 &  20 &  -194 &  52 &  -138 &  -44 &  224 &  8 &  -80 \\ 
& & & & & &  -138 &  -26 &  -132 &  56 &  56 &  68 &  122 &  98 &  -72 &  64 &  148 &  240 \\ 

42. & $a_{1} a_{4} d_{14} e_{6} $ & $\frac{3628}{1815}$ & 25 & 1 &  $\frac{1}{165}$ &  49 &  109 &  69 &  129 &  15 &  -245 &  211 &  -149 &  -165 &  55 &  11 &  119 \\ 
& & & & & &  -231 &  41 &  -79 &  73 &  153 &  37 &  73 &  109 &  -179 &  125 &  109 &  245 \\ 

43. & $a_{1}^2 a_{9} d_{14} $ & $\frac{3608}{1805}$ & 25 & 1 &  $\frac{1}{95}$ &  0 &  76 &  36 &  68 &  0 &  -156 &  112 &  -88 &  -92 &  40 &  8 &  64 \\ 
& & & & & &  -132 &  16 &  -52 &  44 &  84 &  20 &  64 &  60 &  -96 &  72 &  48 &  140 \\ 

44. & $a_{1}^2  a_{9}  d_{14}$ & $\frac{3608}{1805}$ & 25 & 1 &  $\frac{1}{95}$ &  81 &  89 &  -7 &  81 &  13 &  -129 &  33 &  -93 &  -29 &  149 &  5 &  -53 \\
\end{tabular}

\bigskip

\textsc{Table 1 (contd.). Shared vertices} 
\end{center}
\end{sidewaystable}

\begin{sidewaystable}
\begin{center}

\begin{tabular}{c|c|c|>{$}c<{$}|>{$}c<{$}|c|>{$}r<{$}>{$}r<{$}>{$}r<{$}>{$}r<{$}|>{$}r<{$}>{$}r<{$}>{$}r<{$}>{$}r<{$}|>{$}r<{$}>{$}r<{$}>{$}r<{$}>{$}r<{$}}

& name & $\|v\|^2$ & $N$ & $g$ & $\alpha$ & \multicolumn{4}{c}{M} & \multicolumn{4}{c}{O} & \multicolumn{4}{c}{G}\\
\hline
61. & $a_{1}  a_{2} a_{8} e_{6} e_{8}$ & $\frac{3040}{1521}$ & 25 & 1 &  $\frac{1}{117}$ &  20 &  80 &  40 &  88 &  4 &  -172 &  144 &  -108 &  -120 &  44 &  12 &  84 \\ 
& & & & & &  -164 &  36 &  -56 &  60 &  116 &  28 &  60 &  80 &  -132 &  80 &  68 &  172 \\ 

62. & $a_{1} a_{2} a_{4} a_{10} e_{8} $ & $\frac{2968}{1485}$ & 25 & 1 &  $\frac{1}{495}$ &  472 &  516 &  -44 &  404 &  44 &  -736 &  188 &  -428 &  -168 &  692 &  -20 &  -188 \\ 
& & & & & &  -472 &  -92 &  -520 &  288 &  188 &  212 &  472 &  308 &  -188 &  148 &  476 &  832 \\ 

63. & $ a_{1} a_{2} a_{9} d_{13}$ & $\frac{1469}{735}$ & 25 & 1 &  $\frac{1}{210}$ &  -14 &  140 &  68 &  132 &  -18 &  -304 &  262 &  -194 &  -190 &  80 &  14 &  106 \\ 
& & & & & &  -292 &  62 &  -98 &  86 &  182 &  22 &  146 &  120 &  -242 &  170 &  128 &  362 \\ 

64. & $a_{1} a_{2} a_{9} d_{5} e_{8}$ & $\frac{1469}{735}$ & 25 & 1 &  $\frac{1}{210}$ &  202 &  216 &  -20 &  170 &  20 &  -310 &  80 &  -182 &  -72 &  296 &  -8 &  -80 \\
\end{tabular}

\bigskip

\textsc{Table 1 (contd.). Shared vertices} 
\end{center}
\end{sidewaystable}

\begin{sidewaystable}
\begin{center}

\begin{tabular}{c|c|c|>{$}c<{$}|>{$}c<{$}|c|>{$}r<{$}>{$}r<{$}>{$}r<{$}>{$}r<{$}|>{$}r<{$}>{$}r<{$}>{$}r<{$}>{$}r<{$}|>{$}r<{$}>{$}r<{$}>{$}r<{$}>{$}r<{$}}

& name & $\|v\|^2$ & $N$ & $g$ & $\alpha$ & \multicolumn{4}{c}{M} & \multicolumn{4}{c}{O} & \multicolumn{4}{c}{G}\\
\hline

81. & $a_{1}  a_{4}  a_{13} e_{7} $ & $\frac{1748}{875}$ & 25 & 1 &  $\frac{1}{175}$ &  3 &  143 &  55 &  131 &  17 &  -247 &  225 &  -163 &  -177 &  53 &  35 &  93 \\ 
& & & & & &  -249 &  37 &  -69 &  73 &  165 &  7 &  93 &  83 &  -177 &  129 &  95 &  305 \\ 

82. & $a_{1} a_{3}  a_{11} d_{10} $ & $\frac{1732}{867}$ & 25 & 1 &  $\frac{1}{51}$ &  46 &  48 &  -4 &  46 &  8 &  -70 &  20 &  -46 &  -16 &  80 &  4 &  -28 \\ 
& & & & & &  -46 &  -10 &  -48 &  20 &  20 &  24 &  46 &  34 &  -24 &  24 &  52 &  84 \\ 

83. & $a_{1} a_{2} a_{15} d_{7} $ & $\frac{1732}{867}$ & 25 & 2 &  $\frac{1}{51}$ &  -1 &  37 &  17 &  35 &  -5 &  -73 &  63 &  -53 &  -49 &  23 &  7 &  23 \\ 
& & & & & &  -67 &  11 &  -17 &  17 &  47 &  9 &  35 &  25 &  -55 &  41 &  27 &  91 \\ 

84. & $a_{1} a_{9} d_{8} e_{7} $ & $\frac{1688}{845}$ & 25 & 1 &  $\frac{1}{65}$ &  63 &  59 &  -5 &  55 &  5 &  -95 &  25 &  -63 &  -23 &  99 &  3 &  -35 \\
\end{tabular}

\bigskip

\textsc{Table 1 (contd.). Shared vertices} 
\end{center}
\end{sidewaystable}

\begin{sidewaystable}
\begin{center}

\begin{tabular}{c|c|c|>{$}c<{$}|>{$}c<{$}|c|>{$}r<{$}>{$}r<{$}>{$}r<{$}>{$}r<{$}|>{$}r<{$}>{$}r<{$}>{$}r<{$}>{$}r<{$}|>{$}r<{$}>{$}r<{$}>{$}r<{$}>{$}r<{$}}

& name & $\|v\|^2$ & $N$ & $g$ & $\alpha$ & \multicolumn{4}{c}{M} & \multicolumn{4}{c}{O} & \multicolumn{4}{c}{G}\\
\hline
101. & $a_{1} a_{6} a_{9} d_{9}$ & $\frac{629}{315}$ & 25 & 1 &  $\frac{1}{210}$ &  196 &  198 &  -10 &  180 &  30 &  -300 &  70 &  -196 &  -66 &  318 &  6 &  -110 \\ 
& & & & & &  -196 &  -44 &  -210 &  90 &  90 &  86 &  196 &  132 &  -90 &  82 &  210 &  346 \\ 

102. & $a_{1} a_{11}  d_{7} e_{6} $ & $\frac{1248}{625}$ & 25 & 2 &  $\frac{1}{25}$ &  6 &  18 &  10 &  18 &  2 &  -38 &  30 &  -22 &  -26 &  10 &  2 &  18 \\ 
& & & & & &  -34 &  6 &  -14 &  14 &  22 &  6 &  14 &  18 &  -26 &  18 &  14 &  38 \\ 

103. & $a_{1} e_{6}^4$ & $\frac{1248}{625}$ & 25 & 48 &  $\frac{1}{25}$ &  6 &  22 &  10 &  22 &  2 &  -38 &  30 &  -22 &  -22 &  10 &  6 &  18 \\ 
& & & & & &  -34 &  6 &  -14 &  14 &  22 &  2 &  10 &  18 &  -22 &  22 &  14 &  38 \\ 

104. & $a_{1} a_{2}  a_{4}^2 a_{14}$ & $\frac{1208}{605}$ & 25 & 2 &  $\frac{1}{55}$ &  0 &  40 &  16 &  36 &  -4 &  -80 &  72 &  -56 &  -48 &  20 &  8 &  24 \\
\end{tabular}

\bigskip

\textsc{Table 1 (contd.). Shared vertices} 
\end{center}
\end{sidewaystable}

\begin{sidewaystable}
\begin{center}

\begin{tabular}{c|c|c|>{$}c<{$}|>{$}c<{$}|c|>{$}r<{$}>{$}r<{$}>{$}r<{$}>{$}r<{$}|>{$}r<{$}>{$}r<{$}>{$}r<{$}>{$}r<{$}|>{$}r<{$}>{$}r<{$}>{$}r<{$}>{$}r<{$}}

& name & $\|v\|^2$ & $N$ & $g$ & $\alpha$ & \multicolumn{4}{c}{M} & \multicolumn{4}{c}{O} & \multicolumn{4}{c}{G}\\
\hline

121. & $a_{1} a_{9} a_{11} d_{4} $ & $\frac{479}{240}$ & 25 & 1 &  $\frac{1}{120}$ &  104 &  114 &  -8 &  116 &  22 &  -164 &  52 &  -104 &  -38 &  184 &  0 &  -60 \\ 
& & & & & &  -104 &  -26 &  -128 &  60 &  60 &  40 &  104 &  74 &  -46 &  46 &  128 &  198 \\ 

122. & $a_{1}  a_{9}^2 d_{6} $ & $\frac{880}{441}$ & 25 & 4 &  $\frac{1}{21}$ &  2 &  14 &  6 &  14 &  -2 &  -30 &  26 &  -22 &  -22 &  10 &  6 &  10 \\ 
& & & & & &  -26 &  6 &  -6 &  6 &  18 &  2 &  14 &  10 &  -26 &  14 &  10 &  38 \\ 

123. & $a_{1} a_{2}^2  a_{8} e_{6}^2$ & $\frac{880}{441}$ & 25 & 2 &  $\frac{1}{63}$ &  -12 &  48 &  32 &  40 &  -8 &  -88 &  80 &  -60 &  -52 &  24 &  4 &  36 \\ 
& & & & & &  -96 &  16 &  -28 &  16 &  56 &  12 &  44 &  28 &  -60 &  44 &  36 &  112 \\ 

124. & $a_{1} d_{6}^4$ & $\frac{880}{441}$ & 25 & 24 &  $\frac{1}{21}$ &  18 &  18 &  -2 &  22 &  6 &  -26 &  10 &  -22 &  -6 &  30 &  -2 &  -14 \\
\end{tabular}

\bigskip

\textsc{Table 1 (contd.). Shared vertices} 
\end{center}
\end{sidewaystable}

\begin{sidewaystable}
\begin{center}

\begin{tabular}{c|c|c|>{$}c<{$}|>{$}c<{$}|c|>{$}r<{$}>{$}r<{$}>{$}r<{$}>{$}r<{$}|>{$}r<{$}>{$}r<{$}>{$}r<{$}>{$}r<{$}|>{$}r<{$}>{$}r<{$}>{$}r<{$}>{$}r<{$}}

& name & $\|v\|^2$ & $N$ & $g$ & $\alpha$ & \multicolumn{4}{c}{M} & \multicolumn{4}{c}{O} & \multicolumn{4}{c}{G}\\
\hline

141. & $a_{1} a_{2}  a_{4}^2 a_{9}  d_{5}$ & $\frac{269}{135}$ & 25 & 2 &  $\frac{1}{90}$ &  7 &  77 &  33 &  71 &  3 &  -147 &  103 &  -79 &  -91 &  39 &  15 &  65 \\ 
& & & & & &  -129 &  21 &  -45 &  45 &  93 &  15 &  45 &  51 &  -85 &  57 &  45 &  131 \\ 

142. & $a_{1}^4 a_{7}^3$ & $\frac{255}{128}$ & 25 & 6 &  $\frac{1}{32}$ &  26 &  30 &  0 &  30 &  8 &  -44 &  10 &  -32 &  -8 &  46 &  4 &  -20 \\ 
& & & & & &  -28 &  -6 &  -36 &  12 &  12 &  12 &  32 &  18 &  -12 &  12 &  32 &  54 \\ 

143. & $a_{1}^2  a_{4} a_{7}^2 d_{5} $ & $\frac{249}{125}$ & 25 & 2 &  $\frac{1}{50}$ &  46 &  44 &  -6 &  46 &  18 &  -72 &  16 &  -48 &  -14 &  72 &  6 &  -30 \\ 
& & & & & &  -48 &  -6 &  -54 &  18 &  18 &  18 &  48 &  30 &  -18 &  18 &  46 &  84 \\ 

144. & $a_{1} a_{2} a_{3} a_{7}^2 d_{5}$ & $\frac{484}{243}$ & 25 & 4 &  $\frac{1}{27}$ &  -2 &  26 &  14 &  20 &  -2 &  -40 &  34 &  -26 &  -26 &  8 &  6 &  18 \\
\end{tabular}

\bigskip

\textsc{Table 1 (contd.). Shared vertices} 
\end{center}
\end{sidewaystable}

\begin{sidewaystable}
\begin{center}

\begin{tabular}{c|c|c|>{$}c<{$}|>{$}c<{$}|c|>{$}r<{$}>{$}r<{$}>{$}r<{$}>{$}r<{$}|>{$}r<{$}>{$}r<{$}>{$}r<{$}>{$}r<{$}|>{$}r<{$}>{$}r<{$}>{$}r<{$}>{$}r<{$}}

& name & $\|v\|^2$ & $N$ & $g$ & $\alpha$ & \multicolumn{4}{c}{M} & \multicolumn{4}{c}{O} & \multicolumn{4}{c}{G}\\
\hline

161. & $a_{1}^{25}$ & $\frac{31}{16}$ & 25 & 40320 &  $\frac{1}{8}$ &  2 &  8 &  2 &  8 &  0 &  -12 &  8 &  -4 &  -8 &  -2 &  2 &  8 \\ 
& & & & & &  -8 &  0 &  -6 &  2 &  6 &  0 &  2 &  8 &  -8 &  4 &  8 &  12 \\ 

162. & $a_{1}^{25}$ & $\frac{52}{27}$ & 25 & 443520 &  $\frac{1}{9}$ &  3 &  9 &  3 &  9 &  -1 &  -13 &  9 &  -3 &  -9 &  -1 &  1 &  9 \\ 
& & & & & &  -9 &  1 &  -9 &  1 &  7 &  1 &  3 &  9 &  -9 &  3 &  9 &  13 \\ 

163. & $a_{1}^{25}$ & $\frac{48}{25}$
& 25 & 10200960 &  $\frac{1}{5}$ &  1 &  5 &  1 &  5 &  -1 &  -7 &  5 &  -1 &  -5 &  -1 &  1 &  5 \\ 
& & & & & &  -5 &  1 &  -5 &  1 &  5 &  1 &  1 &  5 &  -5 &  1 &  5 &  7 \\ 

164. & $a_{1}^{25}$ & $\frac{48}{25}$
& 25 & 244823040 &  $\frac{1}{5}$ &  4 &  4 &  4 &  4 &  0 &  -8 &  4 &  -4 &  -4 &  8 &  0 &  -4 \\ 
& & & & & &  -4 &  0 &  -4 &  0 &  0 &  0 &  4 &  4 &  -4 &  4 &  4 &  4 \\ 

\end{tabular}

\bigskip

\textsc{Table 1 (contd.). Shared vertices} 
\end{center}
\end{sidewaystable}

\begin{sidewaystable}
\begin{center}

\begin{tabular}{c|c|c|>{$}c<{$}|>{$}c<{$}|c|>{$}r<{$}>{$}r<{$}>{$}r<{$}>{$}r<{$}|>{$}r<{$}>{$}r<{$}>{$}r<{$}>{$}r<{$}|>{$}r<{$}>{$}r<{$}>{$}r<{$}>{$}r<{$}}

& name & $\|v\|^2$ & $N$ & $g$ & $\alpha$ & \multicolumn{4}{c}{M} & \multicolumn{4}{c}{O} & \multicolumn{4}{c}{G}\\
\hline
 
1. & exceptional & $\frac{8}{3}$ & 552 & 495766656000 &  $\frac{1}{3}$ &  2 &  2 &  2 &  2 &  -2 &  -6 &  2 &  -2 &  -2 &  6 &  2 &  2 \\ 
& & & & & &  -2 &  2 &  -2 &  2 &  2 &  2 &  2 &  2 &  -2 &  2 &  2 &  6 \\ 

2. & $\mathrm{HS}_{100}$ $a_1$ & $\frac{9}{4}$ & 101 & 44352000 &  $\frac{1}{4}$ &  1 &  3 &  3 &  3 &  -1 &  -7 &  3 &  -1 &  -3 &  5 &  1 &  3 \\ 
& & & & & &  -5 &  3 &  -3 &  3 &  3 &  1 &  3 &  3 &  -3 &  3 &  3 &  7 \\ 

3. & $\mathrm{HS}_{50}$ $a_2$ & $\frac{32}{15}$ & 52 & 126000 &  $\frac{1}{15}$ &  2 &  10 &  2 &  14 &  -2 &  -22 &  18 &  -14 &  -10 &  2 &  -2 &  14 \\ 
& & & & & &  -22 &  2 &  -2 &  2 &  14 &  6 &  14 &  10 &  -22 &  14 &  6 &  22 \\ 

4. & $J(7,4)$ $a_3$& $\frac{25}{12}$ & 38 & 2520 &  $\frac{1}{12}$ &  7 &  15 &  3 &  15 &  3 &  -17 &  3 &  -13 &  -3 &  15 &  3 &  -11 \\
\end{tabular}

\bigskip

\textsc{Table 2. Additional vertices} 
\end{center}
\end{sidewaystable}

\begin{sidewaystable}
\begin{center}

\begin{tabular}{c|c|c|>{$}c<{$}|>{$}c<{$}|c|>{$}r<{$}>{$}r<{$}>{$}r<{$}>{$}r<{$}|>{$}r<{$}>{$}r<{$}>{$}r<{$}>{$}r<{$}|>{$}r<{$}>{$}r<{$}>{$}r<{$}>{$}r<{$}}

& name & $\|v\|^2$ & $N$ & $g$ & $\alpha$ & \multicolumn{4}{c}{M} & \multicolumn{4}{c}{O} & \multicolumn{4}{c}{G}\\
\hline

21. & $T^3_04^1_0$ $a_{10}$ & $\frac{288}{143}$ & 24 & 1 &  $\frac{1}{143}$ &  128 &  128 &  0 &  120 &  0 &  -216 &  48 &  -128 &  -56 &  216 &  8 &  -48 \\ 
& & & & & &  -128 &  -16 &  -144 &  96 &  48 &  56 &  128 &  88 &  -48 &  64 &  144 &  248 \\

22. & $G_{24,27}$ $a_2^2$ & $\frac{296}{147}$ & 28 & 72 &  $\frac{1}{21}$ &  -2 &  14 &  6 &  14 &  -2 &  -30 &  26 &  -18 &  -18 &  6 &  2 &  10 \\ 
& & & & & &  -30 &  6 &  -6 &  6 &  18 &  2 &  18 &  10 &  -26 &  18 &  10 &  38 \\ 

23. & $G_{22,22}$ $a_1a_4$ & $\frac{161}{80}$ & 27 & 12 &  $\frac{1}{40}$ &  33 &  43 &  -5 &  33 &  1 &  -57 &  17 &  -37 &  -11 &  61 &  -1 &  -13 \\ 
& & & & & &  -37 &  -5 &  -37 &  23 &  13 &  21 &  33 &  29 &  -17 &  15 &  41 &  69 \\ 

24. & $T^{13}_04^1_0$ & $\frac{13864}{6889}$ & 24 & 1 &  $\frac{1}{83}$ &  76 &  76 &  0 &  68 &  0 &  -124 &  28 &  -76 &  -32 &  124 &  4 &  -28 \\
\end{tabular}

\bigskip

\textsc{Table 2 (contd.). Additional vertices} 
\end{center}
\end{sidewaystable}

\begin{sidewaystable}
\begin{center}

\begin{tabular}{c|c|c|>{$}c<{$}|>{$}c<{$}|c|>{$}r<{$}>{$}r<{$}>{$}r<{$}>{$}r<{$}|>{$}r<{$}>{$}r<{$}>{$}r<{$}>{$}r<{$}|>{$}r<{$}>{$}r<{$}>{$}r<{$}>{$}r<{$}}

& name & $\|v\|^2$ & $N$ & $g$ & $\alpha$ & \multicolumn{4}{c}{M} & \multicolumn{4}{c}{O} & \multicolumn{4}{c}{G}\\
\hline

41. & $T^3_20$ $a_{15}$ & $\frac{289}{144}$ & 24 & 1 &  $\frac{1}{48}$ &  43 &  43 &  1 &  41 &  1 &  -73 &  15 &  -43 &  -19 &  71 &  3 &  -19 \\ 
& & & & & &  -43 &  -7 &  -49 &  31 &  17 &  17 &  43 &  27 &  -17 &  21 &  49 &  83 \\ 

42. & $T^8_04^6_0$ & $\frac{289}{144}$ & 24 & 1 &  $\frac{1}{24}$ &  4 &  16 &  8 &  20 &  -2 &  -36 &  26 &  -20 &  -22 &  14 &  2 &  18 \\ 
& & & & & &  -34 &  8 &  -14 &  12 &  20 &  6 &  18 &  16 &  -26 &  18 &  16 &  36 \\ 

43. & $T^2_07^2_0$ $d_8$ & $\frac{580}{289}$ & 25 & 2 &  $\frac{1}{17}$ &  16 &  14 &  -2 &  16 &  2 &  -24 &  10 &  -16 &  -6 &  26 &  2 &  -10 \\ 
& & & & & &  -16 &  -4 &  -18 &  6 &  6 &  6 &  16 &  12 &  -6 &  6 &  18 &  26 \\ 

44. & $T^7_04^7_0$ & $\frac{4432}{2209}$ & 24 & 2 &  $\frac{1}{47}$ &  4 &  36 &  12 &  28 &  0 &  -72 &  60 &  -44 &  -48 &  16 &  0 &  32 \\
\end{tabular}

\bigskip

\textsc{Table 2 (contd.). Additional vertices} 
\end{center}
\end{sidewaystable}

\begin{sidewaystable}
\begin{center}

\begin{tabular}{c|c|c|>{$}c<{$}|>{$}c<{$}|c|>{$}r<{$}>{$}r<{$}>{$}r<{$}>{$}r<{$}|>{$}r<{$}>{$}r<{$}>{$}r<{$}>{$}r<{$}|>{$}r<{$}>{$}r<{$}>{$}r<{$}>{$}r<{$}}

& name & $\|v\|^2$ & $N$ & $g$ & $\alpha$ & \multicolumn{4}{c}{M} & \multicolumn{4}{c}{O} & \multicolumn{4}{c}{G}\\
\hline

61. & $T^{18}_10$ $a_1$ & $\frac{15856}{7921}$ & 24 & 1 &  $\frac{1}{89}$ &  -1 &  61 &  29 &  55 &  -3 &  -129 &  109 &  -89 &  -77 &  37 &  9 &  47 \\ 
& & & & & &  -121 &  25 &  -41 &  33 &  81 &  15 &  65 &  47 &  -105 &  73 &  53 &  151 \\ 

62. & $T^{18}_10$ $a_1$ & $\frac{15856}{7921}$ & 24 & 1 &  $\frac{1}{89}$ &  78 &  90 &  2 &  70 &  14 &  -134 &  30 &  -78 &  -26 &  134 &  -6 &  -46 \\ 
& & & & & &  -78 &  -30 &  -94 &  50 &  46 &  22 &  78 &  38 &  -30 &  34 &  94 &  150 \\ 

63. & $T^1_011^1_0$ $a_6$ & $\frac{2368}{1183}$ & 25 & 2 &  $\frac{1}{91}$ &  82 &  78 &  -2 &  90 &  22 &  -122 &  46 &  -82 &  -26 &  130 &  10 &  -58 \\ 
& & & & & &  -74 &  -18 &  -106 &  38 &  38 &  38 &  82 &  58 &  -38 &  38 &  98 &  150 \\ 

64. & $T^1_011^1_0$ $a_2a_4$ & $\frac{2432}{1215}$ & 25 & 2 &  $\frac{1}{135}$ &  -8 &  92 &  44 &  84 &  -12 &  -196 &  172 &  -128 &  -112 &  44 &  8 &  64 \\
\end{tabular}

\bigskip

\textsc{Table 2 (contd.). Additional vertices} 
\end{center}
\end{sidewaystable}

\end{document}